\documentclass[12pt]{article}
\usepackage{amsmath,amsfonts,amssymb,color,a4,epsfig,graphics}
\usepackage[latin1]{inputenc}

\newcommand{\R}{{\mathbb{R}}}
\newcommand{\C}{{\mathbb{C}}}
\newcommand{\Z}{{\mathbb{Z}}}
\newcommand{\N}{{\mathbb{N}}}

\def\pa{\partial}
\def\ra{\rightarrow}
\def\preuve{\begin{proof}}
\def\ga{\alpha}
\def\gb{\beta}

\def\gd{\delta}
\def\ge{\varepsilon}

\def\gg{\gamma}

\def\gl{\lambda}
\def\go{\omega}

\def\gs{\sigma}

\def\gO{\Omega}

\newtheorem{defi}{Definition}[section]

\newtheorem{lemm}{Lemma}[section]
\newtheorem{prop}{Proposition}[section]
\newtheorem{rem}{Remark}[section]

\newtheorem{theo}{Theorem}[section]
\newtheorem{exem}{Example}[section]
\newenvironment{demo}{\noindent {\it Proof.--}
      \begin{quotation}\noindent}{\end{quotation}\hfill$\square $}

\renewcommand{\footnote}[1]{\footnotetext{#1}}

\begin{document}

\title{Essential self-adjointness for combinatorial Schr\"odinger
operators \\ III- Magnetic fields}
\author{Yves Colin de Verdi\`ere \footnote{Grenoble University,
Institut Fourier,
 Unit{\'e} mixte
 de recherche CNRS-UJF 5582,
 BP 74, 38402-Saint Martin d'H\`eres Cedex (France);
{\tt yves.colin-de-verdiere@ujf-grenoble.fr};
{\tt http://www-fourier.ujf-grenoble.fr/$\sim $ycolver/}
}\\ Nabila Torki-Hamza \footnote{Universit\'e du 7  Novembre \`a Carthage,
Facult\'e des Sciences de Bizerte,
Math\'ematiques et Applications (05/UR/15-02),
 7021-Bizerte (Tunisie);
{\tt nabila.torki-hamza@fsb.rnu.tn};
{\tt torki@fourier.ujf-grenoble.fr}}
\\ Fran{\c c}oise Truc \footnote{Grenoble University, Institut Fourier,
Unit{\'e} mixte
 de recherche CNRS-UJF 5582,
 BP 74, 38402-Saint Martin d'H\`eres Cedex (France);
{\tt francoise.truc@ujf-grenoble.fr};
{\tt http://www-fourier.ujf-grenoble.fr/$\sim$trucfr/}
}}

\maketitle
\begin{abstract}

We define the magnetic Schr\"odinger operator on an infinite  graph
 by the data of a magnetic field, some  weights on
vertices and some  weights on edges. We discuss essential self-adjointness
of this operator for  graphs  of bounded degree.
The main result is a discrete version of a result of
two authors of the present paper. \\

\emph{On d\'efinit l'op\'erateur de Schr\"odinger avec champ magn\'etique
 sur un graphe
infini  par la donn\'ee d'un champ magn\'etique, de poids sur les sommets
et de poids sur les ar\^etes. Lorsque le graphe est de degr\'e born\'e,
on \'etudie le caract\`ere essentiellement auto-adjoint d'un tel
op\'erateur.
Le r\'esultat principal est une version discr\`ete d'un r\'esultat de
deux  des auteurs du pr\'esent article.\\}


\textbf{Math Subject Classification (2000):} 05C63, 05C50, 05C12, 05C05, 35J10, 47B25.
\end{abstract}

\newpage

\section{Introduction}
\paragraph{}
In this work, we investigate essential self-adjointness
for magnetic Schr\"odinger operators on infinite
weighted graphs $G=(V,E)$ of bounded degree.
It is a continuation of \cite{To} and \cite{ColToTr}, where the
 same problem was studied in the non
magnetic case, for metrically complete graphs (\cite{To})
as well as non complete ones (\cite{ColToTr}). The main result is a
discrete version of the result in \cite{ColTr}.

In the former paper (\cite{ColToTr}), we proved that provided a
growth condition on the potential $W$, namely $W\geq N/(2D^2)$
where $D$ is the distance to infinity and $N$ the maximal degree,
 the  Schr\"odinger operator
$\Delta_{\go,c}+ W$ is essentially self-adjoint.
The operator $\Delta_{\go,c}$ is defined, for any weights
$\go:V\longrightarrow \R_+^\star$ and $c:E\longrightarrow \R_+^\star,$
 by:
$$ \left( \Delta_{\go,c}f\right)\left( x\right)
 =\dfrac{1}{\go_{x}^{2}} ~\sum_{y\sim x}
 ~c_{\lbrace x,y\rbrace}\left( f\left( x\right) -f\left(y \right) \right),  $$
for any $~f\in C_0(V)~$ (finite supported function) and any vertex $x\in V .$

We will extend  this result
  to the case of magnetic graph Laplacians. To
do this, we establish a lower bound for the magnetic
 Dirichlet integral, in terms of
an effective potential depending on the magnetic field
 (and not depending on the magnetic potential) and
follow the method  described in \cite{ColToTr}. The precise
setup of the result is
described in the next section.

\section{Magnetic fields on graphs}
\paragraph{}
Discrete magnetic Schrödinger operators were already introduced by
several authors, see \cite{Li-Lo,Col,Col3,Tor}.

\subsection{Magnetic Schr\"odinger operators}
\paragraph{}
Let $G=(V,E) $ be  a {\it locally finite  connected  graph}.
We will denote by $\{ x,y\}\in E  $ an edge
and by $[x,y]$ and $[y,x]$ the two orientations of this edge.
\newline
We equip $G$  with
\begin{description}
  \item[(i)] a  set  of non zero {\it complex}
weights on oriented edges:
 $C_{xy}\in \C \setminus 0$ for $ \{x,y \} \in E $ with
$C_{yx}=\overline{C_{xy}}$. We write
$ C_{xy}= c_{xy}e^{i \ga _{xy}}$ with $c_{xy}=|C_{xy}|.$
We have $c_{xy}=c_{yx}$ and
 $ \ga _{xy}=- \ga _{yx}.$
The set  $A =(\ga _{xy})$ is called a {\it magnetic potential} on $G.$
  \item[(ii)] a  set of {\it strictly positive}
weights on the  vertices:  $\go _x ,~x\in V.$
\end{description}

 The space of complex valued functions on the graph $G$ is denoted here by
$$C(V)=\lbrace f:V \longrightarrow {\C}\rbrace$$  and
 $C_0(V)$ is the subspace of $C(V)$ of functions with finite support.\\

 We consider the Hilbert space
  $$l^2_\go (V)=\lbrace f\in C(V);\; \underset{x\in
V}{\sum}\go_x^2~| f(x)|^2 <\infty \rbrace$$
equipped with the Hermitian inner product given by
 \[ \langle f,g\rangle _{l^2_\go}=\sum _{x\in V}\go_x^2f(x)\overline{g(x)} .\]
Let us consider the Hermitian form
\[ Q_{c,A} (f)=\sum _{\{x,y \} \in E}c_{xy} | f(x) -e^{i \ga
  _{xy}}f(y)|^2 ,\]
where we take only one term for each (unoriented) edge (the contribution
is the same for both choices of orientations $[x,y]$ and $[y,x$].)

The associated \emph{magnetic Schr\"odinger operator}
$H_{\go,{c,A}} $ is given formally by
\[ \langle H_{\go,{c,A}}f|f \rangle _{l^2_\go}=
Q_{c,A} (f) .\]
We get easily
\[ H_{\go,{c,A} }f(x)=\frac{1}{\go_x^2}  \sum _{y\sim x }
c_{xy}[ f(x)-e^{i \ga
  _{xy}}f(y)].\]

This operator $H_{\go,{c,A}} $ is Hermitian symmetric on $C_0(V)$
with the Hermitian product induced by $l^2_\go (V).$

\subsection{Gauge transforms}
\begin{defi}
Let us consider a sequence of complex numbers
 $(u_x)_{x\in V}$ with $|u_x|\equiv 1 $
and write $u_x=e^{i\sigma_x}$.
The associated \emph{gauge transform} $U$ is the unitary map on
 $l^2_\go  (V)$
defined by
$$(Uf)(x)= u_x f(x) .$$
The map $U $ acts on the quadratic forms  $Q_{c,A}$
by $$U^\star (Q_{c,A})(f) = Q_{c,A} (Uf).$$
Let us define the magnetic potential $U^\star (A) $ by
$U^\star (Q_{c,A})=Q_{c,U^\star (A)}\ .$
  The associated magnetic Schr\"odinger operator
is $ H_{\go ,c,U^\star (A)}$.
\end{defi}
More explicitly, we get:
\[U^\star (A)_{xy}=  \ga _{xy}+\sigma_y
  -\sigma_x. \]
Let us denote by $C_1 (G)$ the $\Z$-module generated by
the oriented edges with the relation $[x,y]=-[y,x]$,
and by $\pa$ the boundary operator
$$C_1 (G)\overset{\pa}{\ra}C(V,\Z) $$ so that
$\pa ([x,y])=\gd _y -\gd _x $
where, for $z\in V$, $\gd_z  \in C(V,\Z) $
is defined by $\gd_z(z)=1$ and $\gd _z(z')=0$ if $z'\ne z.$
We will denote, for any $\gg \in C_1 (G),~\gg =\sum_{e\in E}\gg (e)e.$\\

The space of cycles, denoted by $Z_1(G)$,  is the kernel of the boundary
operator. If $G$ is connected and finite, it is a known fact
 (see \cite{Bi}, part 1, chapters 4, 5)
that $Z_1(G)$ is a free $\Z-$module,
of rank $\# V -\# E +1$,  with a basis of geometric cycles
 $\gg =[x_0, x_1]+[x_1,x_2]+ \cdots  + [x_{n-1} ,x_n ]$
  with, for $i=0, \cdots, n-1$, $\{ x_i ,x_{i+1} \} \in E ,$
  and $x_{n}=x_{0}.$\\

  We will construct a basis of cycles for any graph $G.$

  Zorn's Lemma allows to show the existence of  a maximal tree $T$
of $G.$
 We have $V(T)=V(G)$ and we  denote $E'$
the set of all edges of $G$ which are not edges of $T$. We choose an
orientation for each edge of $E'$.
 For any  (oriented) edge  $[x,y] \in E'$, there exists a unique simple
 path $\gb_{yx}$ in the tree $T$ linking $y$ to $x.$
 So $\gg_{xy}=[x,y]+\gb_{yx}$ is a geometric
 cycle of $G$.\\
 Let $\gg \in Z_1(G)$, we set:
 $\displaystyle \gg'=\gg - \sum_{e\in \gg \cap E'}\gg (e)\gg _{e},$ where
 we denote $\gg _{e}=\gg_{xy},$ if $e=[x,y]$ is an oriented edge which is
 in $E'$ and included in the cycle $\gg.$
 Then  $\gg'$ is a cycle of $G$ with support in $T.$
 So it vanishes and we have:
 $\displaystyle \gg=\sum_{e\in \gg \cap E'}\gg (e)\gg _{e}.$\\
 We have proved the following Lemma:

\begin{lemm} \label{lemm:basecyc}
Let $T$ a maximal tree of $G.$ To each oriented edge $[x,y]$ of $G,$
if $\{x,y\}$ is not an edge of $T,$ we have  associated a unique cycle
$\gg_{xy}$ of the graph $G$ including $[x,y].$
The set of all such cycles is a basis of $Z_1(G).$
\end{lemm}

\begin{defi}
Let us define the {\rm holonomy map}: $${\rm Hol}_A :Z_1(G)\ra \R/2\pi \Z $$
by
$${\rm Hol}_A ([x_0, x_1]+[x_1,x_2]+ \cdots  + [x_{n-1},x_0 ] )=
\ga _{x_0x_1}+ \cdots +\ga _{x_{n-1} x_0} .$$
\end{defi}
\begin{prop}
With the notations above we have:
\begin{description}
  \item[\emph{(i)}] The map $A \ra {\rm Hol}_A $ is surjective
onto  ${\rm Hom }_\Z(Z_1(G),\R/2 \pi \Z ).$
  \item[\emph{(ii)}] ${\rm Hol}_{A _1} ={\rm Hol}_{A _2}  $
if and only if there exists a gauge transform $U$
so that $U^\star (A _2 )= A _1$.
\end{description}
\end{prop}

\begin{demo}
For (i), let $hol \in {\rm Hom }_\Z(Z_1(G),\R/2 \pi \Z )$ and
let $T$ a maximal tree.
We choose $A=(\ga_{xy})$ such that $\ga_{xy}=hol(\gg_{xy})$
if $\{x,y\}\in E'$ and  $\ga_{xy}=0$ if $\{x,y\}\in E(T),$
see Lemma \ref{lemm:basecyc}.\\
Then we have $hol=Hol_{A}.$ \\
For (ii), it suffices to prove that, if $Hol_{A}=0,$ then there
exists a gauge transform $U $ so that  $U^\star (A )=0.$\\
We must find a sequence $(\gs_{x})$ satisfying the equality
$\gs _{x }=\ga_{xy} +\gs_{y},$ for any edge $\{x,y\}.$\\
We fix $x_{0}\in V$, set $\gs _{x _{0}}=0$. If
$x\in V\setminus \{x_{0}\},$ there exists a path
$x_{0},x_{1},...,x_{n-1},x_{n}=x$ connecting $x_{0}$ to $x$, and we set
$\gs_{x}=\ga_{xx_{n-1}}+....+\ga_{x_{2}x_{1}}+\ga_{x_{1}x_{0}}.$
This doesn't depend on the path  from $x_{0}$ to $x$, since $Hol_{A} (\gg)=0$ for any
cycle $\gg.$
\end{demo}

In the case of finite planar graphs, the assertions (i) and (ii)
are similar to respectively Lemma 2.2 and Lemma 2.1 in \cite{Li-Lo}.

\begin{defi}
A {\rm magnetic field} $B$ on the graph $G$ is given by
an holonomy map $${\rm hol} \in{\rm Hom }_\Z(Z_1(G),\R/2 \pi \Z ).$$
If $B$ is associated to the {\rm magnetic potential} $A$, we write
$B=dA$ and we have ${\rm hol}={\rm Hol}_A$.
\end{defi}
\begin{rem}
The magnetic Schr\"odinger operator $H_{\go,c,A} $
is uniquely defined, up to unitary conjugation,
by the data of the magnetic field $B$, the weights $(c_{xy})_{\{x,y\}
\in E}$ and the weights $\go =(\go_x)_{x\in V}$.
\end{rem}
\subsection{Norms of magnetic fields}

\begin{defi}
If $G=(V,E)$ is a finite connected graph with a magnetic field
$B$, we define the norm $|B|$ of $B$ as the lowest eigenvalue of $H_{1,1,A}$
on $l^2_{\go} (V)$ with $\go\equiv 1$, for any $A$ with $dA =B$.
\end{defi}

\begin{lemm}
We have $| B|=0$ if and only if ${\rm Hol}_A =0.$
\end{lemm}
\begin{demo}
If ${\rm Hol}_A =0$,  $H_{1,1,A}$ is unitarily equivalent
to $H_{1,1,0}$ whose lowest eigenvalue is $0$ with constant
eigenfunctions.\\
Conversely, let $f\ne 0 $ with $ H_{1,1,A}f=0$
and hence $Q_{1,A}(f)=0.$
This implies that all terms in the expression of
$Q_{1,A}(f)$
vanish~: for any edge $\{x,y\}$ we have $f(x) =e^{i \ga
 _{xy}}f(y)$. If $\gg=[x_0, x_1]+[x_1,x_2]+ \cdots  + [x_{n-1},x_0 ]$ is a cycle,
we have in particular
\[  f(x_n)=e^{i\ga_{x_{n}x_{n-1}}}f(x_{n-1})=
\cdots =e^{-i {\rm  Hol}_A (\gg)} f(x_0)~.\]
Hence
$$e^{-i {\rm  Hol}_A (\gg)}=1.$$
\end{demo}
\begin{lemm}\label{lemm:cyclic}
Let $G=\Z/N\Z$ be the cyclic graph with $N$ vertices,
 $\gO$ be  the holonomy of
 $\gamma=[0,1]+ [1,2]+ \cdots [N-1,0]$
and $\gd =\inf _{k\in \Z} |\gO -2\pi k|$.\\
If $B$ denotes the magnetic field such that $B(\gamma) = \gO$ we have
\[ |B| =  |1-e^{i\gd/N}|^2.\]
In particular, the maximal value
\[ |B| = |1-e^{i\pi/N}|^2 \]
is obtained for $B =\pi $.
\end{lemm}
\begin{demo}
We can choose $A$ so that
$$Q_{1,A}(f)=\sum _{x=0}^{N-1} |f(x)-e^{i\Omega /N} f(x+1) |^2 .$$
The eigenvectors
are the $N$ complex functions on $V$:
$$f_\xi: x\ra \xi^x $$
where $\xi ^N=1.$
We have $\| f_\xi \| _{l^2}^2 =N $
and
$$Q_{1,A}(f_\xi)= N |1 -\xi e^{i\Omega /N}|^2.$$
\end{demo}

\subsection{Lower bounds using an effective potential}
\begin{defi}
Let $m\in \N$.
A {\rm good covering of degree $m$} of $G=(V,E)$
is a family $G_l =(V_l, E_l )$ with $l \in L$
 of finite connected sub-graphs of
$G$ so that
\begin{description}
  \item[\emph{(i)}]$ V=\cup_{l\in L} V_l $
  \item[\emph{(ii)}] for any $\{x, y\} \in E,$
\[ 0< \# \{ l \in L~|~ \{ x, y \} \in E_l  \} \leq m .\]
\end{description}
\end{defi}
\begin{exem}
Let $G$ the 1-skeleton of a triangulation of the plane $\R^{2}.$
Then the set of all the triangles of this triangulation is a
good covering of degree $2$.
\end{exem}
\begin{rem}
A graph $G$  of bounded degree admits   good coverings
  given by the:

\begin{prop}
Let $G=(V,E)$ be a graph of bounded degree $N$. For $k\geq 1$
and $x\in V$,  let
$$G_{x}^{k}=\{y\in V,~d(x,y)\leq k\}.$$
The family $(G_{x}^{k})_{x\in V}$ is a good covering of $G$ of
degree $m=\frac{N(N-1)^{k}-2}{N-2}$ of the graph $G.$
\end{prop}
\end{rem}
The main estimate is given by the following Theorem:
\begin{theo}\label{theo:ESA}
Let $(G_l)_{l \in L} $ a good covering of degree $m$  of  $G.$
Then for any $f\in C_0(V),$
$$ Q_{c,A} (f) \geq \sum_{x\in V} W(x) \go_x^2 |f(x)|^2$$
with the {\rm effective potential}
\begin{equation}\label{QW}
 W(x)= \frac{1}{m}  \sum _{\{ l\in L | x\in V_l \}} |B_l |
\inf _{\{y,z\} \in E_l}
c_{yz}
\end{equation}
 where $ |B_l |$ is the norm of the restriction of $B$
to $G_l .$
\end{theo}
\begin{demo}
From the definition of $Q_{c,A}$ and $m,$ we have
\[ Q_{c,A}(f)\geq \frac{1}{m} \sum _{ l \in L }
\sum _{\{x,y\} \in E_l } c_{xy} |f(x)-e^{i \ga_{xy}}f(y)|^2. \]
Using the definition of $ |B|_l ,$ we get
\[ Q_{c,A}(f)\geq \frac{1}{m}\sum _l \left( \inf_{\{y,z\} \in E_l}
c_{yz}\right) |B|_l
\left( \sum _{x\in V_l } \go_x^2  |f(x)|^2
\right)
\]
which gives the lower bound.
\end{demo}

\section{Magnetic confinement}
\paragraph{}
  We want to find a criterion similar to the main result
of \cite{ColTr} which says that if $(c,| B |) $ grows fast enough
near infinity, then $H_{\go,c,B }$ is essentially self-adjoint
on $C_0(V)$.
We will use Theorem 4.3 of \cite{ColToTr} which gives the case where $B=0.$

Let us define the distance
$d_p$ given in terms of the weights $\go$ and $c$ as follows:
\[ p_{xy}= \frac{\min (\go_x, \go_y)}{\sqrt{c_{xy}}} \]
for any vertices $x,y \in V.$

 Define also $D(x)~ (\leq \infty) $ (see \cite{ColToTr})
as the distance from a vertex $x$ to the boundary $V_\infty.$
\begin{theo}\label{neni}
Let $G$ a graph with maximal degree $N$ and let $(G_l)_{l \in L} $ a good
covering of degree $m$ of $G.$
If there exists $M$ so that $$ W (x) \geq {N
}/{2 D(x)^2}-M ,$$ where $W$ is the effective potential
defined in (\ref{QW}), then $H_{\go,c,B}$ is essentially self-adjoint.
\end{theo}
\begin{rem}
 Theorem \ref{neni} holds in particular if $(G,d_p)$ is a complete metric space.
\end{rem}

\begin{demo} The proof follows the steps of the proof
of Theorem 4.3. in \cite{ColToTr}. \\
In particular we use the following Agmon estimate.

\begin{lemm}\label{nen}
Let $v$ be a weak solution of $Hv=0$, and let $f=\overline{f} \in C_{0}(V)\ $
 a real function with finite support in $V.$
Then
 \begin{equation}\label{ute}
\langle   fv\ , H (fv) \rangle  \ =\  \frac{1}{2}\sum_{x\in V}\sum_{y\sim x } {\Re} [v(x) \overline v(y) C_{yx}](f(x)-f(y))^2
\end{equation}

\end{lemm}

\begin{demo}
We denote here $H\equiv H_{\go,c,B}.$ \\
The proof is a simple calculation:
\begin{align}
\langle fv\ , H (fv) \rangle
&=\sum _{x\in V}f(x)v(x)\left( \sum_{y\sim x}\ c_{xy}[f(x)\overline v(x)- e^{-i \ga _{xy}}f(y)\overline v(y)] \right.
\notag \\
&+\left.\sum_{y\sim x}\ W(x) f(x)\overline v(x)\right)
\notag \\
&=\ \sum_{x\in V}f(x)v(x)\left( \sum_{y\sim x}\ \overline {C_{xy}}(f(x)-f(y))\overline v(y)\right)
\notag
\end{align}
where we used the fact that $\overline{H} \overline{ v}=0$.\\
An edge $\{x,y\}$ contributes to the first sum twice.
So the total contribution is
$$f(x)v(x)\ \overline {C_{xy}}(f(x)-f(y))\overline v(y)- f(y)v(y) C_{xy}(f(x)-f(y))\overline v(x)$$
so
$$\langle   fv\ , H(fv) \rangle =\sum_{\{x, y\} \in E} [f(x)-f(y)]\left[f(x)v(x) C_{yx}\overline v(y) -f(y)v(y)  C_{xy}\overline v(x)\right] ~$$
Noticing that the quantity is real, we take the mean value of the expression and of its conjugate then we get the result.
\end{demo}

 From Lemma \ref{nen} we derive the following Theorem.
\begin{theo}\label{Horo}
Let $v$ be a  solution of $(H-\gl)v=0$.
 Assume that $v$ belongs to $l^2_{\go}(V )$
and that there exists a constant $c>0$ such that, for all $u \in C_0(V)$,
\begin{equation}\label{boun}
\langle u | (H-\gl) u \rangle_{l^2_{\go}}  \geq
 \frac{N}{2}\sum_{x\in V}\max \left(\frac{1 }{D(x)^2},1
\right) \go_x^2  |u(x)|^2  +  c \|u\|_{l^2_{\go}}^2,
\end{equation}
then $v\equiv 0$.
\end{theo}
\begin{demo}
We refer to \cite{ColToTr} for the proof, since the fact that we use complex functions does
not make any change in it.
\end{demo}

Then Theorem \ref{neni} follows from Theorems \ref{theo:ESA} and \ref{Horo}
 since we have for any $u \in C_0(V)$ :

$$\langle u | H u \rangle_{l_{\go}^2} \geq \sum_{x\in V} W(x) \go_x^2  |u(x)|^2, $$
so
$$\langle u | (H-\gl) u \rangle_{l^2_{\go}}  -\frac{N}{2}\sum_{x\in V} \frac{1}{D(x)^2}
\go_x^2   |u(x)|^2
\geq   \sum_{x\in V}\
 -(M+\gl) \|u\|_{l^2_{\go}}^2  . $$

\end{demo}

\section{Examples}
The simplest example where we can make estimates is
the ``infinite ladder'', see \cite{Li-Lo}. That is
the graph $G=(V,E)$ where the set of vertices $V$ is the Cartesian product
$V=\N \times \{ -1,1 \} $ equipped with the ``horizontal'' edges
$\{ (l,\ge ), (l+1,\ge ) \}$ with $l=0,1, \cdots$ and $\ge =\pm 1$
and the ``vertical'' edges $\{ (l,-1),(l, +1) \} $
with $l=0,1, \cdots.$
We will use the ``square'' cycles
$$\gg_l= [(l,1),(l+1, 1)]+[(l+1, 1),(l+1, -1)]+[ (l+1,-1),(l,-1)]+[(l,-1),(l,1)],$$
for $l=0,1,\cdots$, as a basis of the space of cycles.
Let $b_l $ be the holonomy of $B$ in the cycle $\gg_l$.
We will take $B$ so that the value of $|B_{\gg_l}|$ is
$2-\sqrt{2}$, which is the maximal one by Lemma \ref{lemm:cyclic}.

Using the good covering of $G$ by the cycles $\gg_l$, we get $m=2$ and
the effective potential
$$W((l,\ge ))= \left(1-\frac{\sqrt{2}}{2}\right)\inf_{ \{ x,y\} \in E_l} c_{xy}.$$
We will take
$$c_{ (l,\ge)(l+1,\ge)}= c_{(l,-1)(l,+1)}=C_l ~
\textnormal{and}~~\go_{(l,\ge) }=w_l~ .$$
If $C_l$ is increasing, we get
$W((l,\ge ))= \left(1-\frac{\sqrt{2}}{2}\right)C_l.$
Let us assume that $w_l $ is decreasing,
we get
$$p_{ (l,\ge)(l+1,\ge)}= \frac{\go_{l+1}}{\sqrt{C_l}}~
\textnormal{and}~~
D( (l,\ge))=\sum_{m=l}^\infty   \frac{w_{m+1}}{\sqrt{C_m}}.$$
We take now $C_l =l^a $ with $a>0$
and $w_l=l^{-b}$ with $b>0 .$ The graph is not complete for
the distance $d_p$ if $a+b/2>1.$
In this case, we have
$$D( (l,\ge))\sim c_1  l^{(1-a-b/2)}~
\textnormal{and}~~ W((l,\ge ))\sim c_2 l^a.$$
The assumption of Theorem \ref{theo:ESA} is satisfied
if
$b<1.$ So we get that, if $0 <b <1 $ and $a+b/2 >1,$ the operator
$H_{\go,c,B}$ is essentially self-adjoint.
If $a>2 ,$ the operator
$H_{\go,c,0}$ is not essentially self-adjoint by Theorem 4.1 in \cite{ColToTr}.

\section{Questions}

The following questions are unsolved at the moment:

\begin{enumerate}
\item If $H_{\go,c,0} $ is essentially self-adjoint, does it imply that
$H_{\go,c,B}$ is essentially self-adjoint for any $B$?
Does it hold in the continuous case?
\item What would be a correct statement for a locally finite
graph with unbounded degree (even if $B=0$)?
\item In the case where the completion of $(G,d_p)$
is compact and under the assumption of Theorem
\ref{theo:ESA},  $H_{\go,c,B}$ has a compact resolvent.
What is the asymptotic behavior of the eigenvalues?
The continuous case is worked out in \cite{Col2}.
\end{enumerate}


\begin{thebibliography}{666}

\bibitem[Bi]{Bi} N.Biggs:
\newblock {\it Algebraic Graph Theory,}
\newblock Cambridge University Press (1974).

\bibitem[B-M-S]{B-M-S} M. Braverman, O. Milatovic \&  M. Shubin:
\newblock Essential self-adjointness of
Schr\"odinger-type operators on manifolds,
\newblock {\it Russian Math. Surveys} {\bf 57}, 641--692  (2002).

\bibitem[CdV]{Col} Y. Colin de Verdi\`ere:
\newblock {\it Spectre de graphes,}
\newblock Cours sp\'ecialis\'es {\bf  4},
\newblock Soci\'et\'e math\'ematique de France  (1998).

\bibitem[CdV2]{Col2} Y. Colin de Verdi\`ere:
\newblock {\it Asymptotique de Weyl pour les bouteilles magn\'etiques,}
\newblock Commun. Math. Phys. {\bf 105}, 327--335  (1986).

\bibitem[CdV3]{Col3} Y. Colin de Verdi\`ere:
{\it Multiplicities of eigenvalues and tree-width of graphs},
    {J. Combin. Theory Ser. B},
  {\bf 74}, {121--146} (1998).


\bibitem[CdV-Tr]{ColTr} Y. Colin de Verdi\`ere \& F. Truc:
Confining quantum particles with a purely magnetic field,
\textcolor{blue}{\tt ArXiv:0903.0803v3} [math-ph],
\newblock {\it Ann. Inst. Fourier (Grenoble)} (to appear)(2010).

\bibitem[CdV-To-Tr]{ColToTr} Y. Colin de Verdi\`ere, N. Torki-Hamza
\& F. Truc:
\newblock {\it Essential self-adjointness for combinatorial
Schr\"odinger operators II. Metrically non complete graphs,}
 \newblock \textcolor{blue}{\tt  ArXiv:1006.5778v2}[math.SP], (2010),
to appear in ``Mathematical Physics, Analysis and Geometry''.

\bibitem[Dod]{Dod} J. Dodziuk:
\newblock Elliptic operators on infinite graphs,
\newblock {\it Analysis geometry and topology of elliptic operators,}
353--368, World Sc. Publ., Hackensack  NJ. (2006).

\bibitem[Du-Sc]{Du-Sc} N. Dunford \&  J. T. Schwartz:
\newblock {\it Linear operator  II, Spectral Theory,}
\newblock  John Wiley {\rm \&} Sons, New York  (1971).

\bibitem[Li-Lo]{Li-Lo}  E. Lieb \& M. Loss:
  \newblock   {Fluxes, {L}aplacians, and {K}asteleyn's theorem},
  \newblock {\it {Duke Math. J.},}
     {\bf 71}, {337--363}  (1993).

\bibitem[Nen-Nen]{Nen} G. Nenciu \& I. Nenciu:
\newblock  On confining potentials and essential self-adjointness for
Schr\"odinger operators on bounded domains in $\R^n$,
\newblock {\it Ann. Henri Poincar\'e,} {\bf 10}, 377--394 (2009).

\bibitem[Ol]{Ol} I.M. Oleinik:
\newblock On the essential self-adjointness
of the Schr$\ddot{o}$dinger operator on complete Riemannian manifolds,
\newblock  {\it Mathematical Notes} {\bf 54}, n 3, 934--939 (1993).

\bibitem[R-S]{RS} M.Reed \&  B.Simon:
\newblock {\it Methods of Modern mathematical Physics I,}
 {\it Functional analysis,} (1980),
\newblock {\it II, Fourier analysis, Self-adjointness} (1975),
\newblock  New York Academic Press.

\bibitem[Shu]{Shu} M. Shubin:
Essential self-adjointness
for semi-bounded magnetic Schr\"odinger operators on
non-compact manifolds,
\newblock {\it J. Func. Anal.,} {\bf 186}, 92--116  (2001).

\bibitem[Sh]{Sh} M. Shubin:
\newblock {\it Classical and quantum completness
 for the Schr\"{o}dinger operators on non-compact manifolds,}
\newblock Geometric Aspects of Partial Differential
\newblock Equations (Proc. Sympos., Roskilde, Denmark (1998))
\newblock Amer. Math. Soc. Providence, RI, 257--269  (1999).

\bibitem[To]{To} N. Torki-Hamza:
\newblock Laplaciens de graphes infinis I- Graphes m\'etriquement complets,
\newblock {\it Confluentes Mathematici}, {\bf 2}, n3, 333--350  (2010).

\bibitem[Tor]{Tor} N. Torki-Hamza:
\newblock {\it Stabilit\'e des valeurs propres avec champ magn\'etique
sur une vari\'et\'e Riemannienne et sur un graphe,}
\newblock Th\`ese de doctorat de l'Universit\'e de Grenoble I, France,
(1989).

\bibitem[Wo]{Wo} R.K. Wojiechowski:
\newblock {\it Stochastic completeness of graphs,}
\newblock Ph.D. Thesis,
\newblock The graduate Center of the University of New-York  (2008).

\end{thebibliography}
\end{document}